\documentclass{article}
\usepackage{harvard}
\usepackage{graphicx}
\usepackage{camnum}
\usepackage{amsmath}

\newcommand{\ee}{{\mathrm e}}
\newcommand{\ii}{{\mathrm i}}

\begin{document}
\title{Why geometric numerical integration?}

\author{A. Iserles \& G.R.W. Quispel}

\maketitle

\thispagestyle{empty}

\section{The purpose of GNI}

Geometric numerical integration (GNI) emerged as a major thread in numerical mathematics some 25 years ago. Although it has had antecedents, in particular the concerted effort of the late Feng Kang and his group in Beijing to design structure-preserving methods, the importance of GNI has been recognised and its scope delineated only in the 1990s. 

But we are racing ahead of ourselves. At the beginning, like always in mathematics, there is the definition and the rationale of GNI. The rationale is that all-too-often mathematicians concerned with differential equations split into three groups that have little in common. Firstly, there are the applied mathematicians, the model builders, who formulate differential equations to describe physical reality. Secondly, there are those pure mathematicians investigating differential equations and unravelling their qualitative features. Finally, the numerical analysts who flesh out the numbers and the graphics on the bones of mathematical formulation. Such groups tended to operate in mostly separate spheres and, in particular, this has been true with regards to computation. Discretisation methods were designed (with huge creativity and insight) to produce rapidly and robustly numerical solutions that can be relied to carry overall small error. Yet, such methods have often carried no guarantee whatsoever to respect qualitative features of the underlying system, the very same features that had been obtained with such effort by pure and applied mathematicians.

Qualitative features come basically in two flavours, the {\em dynamical\/} and the {\em geometric.\/} Dynamical features -- sensitivity with respect to initial conditions and other parameters, as well as the asymptotic behaviour -- have been recognised as important by numerical analysts for a long time, not least because they tend to impinge directly on accuracy. Thus, sensitivity with respect to initial conditions and perturbations comes under `conditioning' and the recovery of correct asymptotics under `stability', both subject to many decades of successful enquiry. Geometric attributes are invariants, constants of the flow. They are often formulated in the language of differential geometry (hence the name!) and mostly come in three varieties: {\em conservation laws,\/} e.g.\ Hamiltonian energy or angular momentum, which geometrically mean that the solution,  rather than  evolving in some large space $\BB{R}^d$,  is restricted to a lower-dimensional manifold $\mathcal{M}$, {\em Lie point symmetries,\/} e.g.\  scaling invariance, which restrict the solution to the tangent bundle of some manifold, and quantities like {\em symplecticity\/} and {\em volume,\/} whose conservation corresponds to an evolution on the cotangent bundle of a manifold. {\em The design and implementation of numerical methods that respect geometric invariants is the business of GNI.\/}

Since its emergence, GNI has become the new paradigm in numerical solution of ODEs, while making significant inroads into numerical PDEs. As often, yesterday's revolutionaries became the new establishment. This is an excellent moment to pause and take stock. Have all the major challenges been achieved, all peaks scaled, leaving just a tidying-up operation? Is there still any point to GNI as a separate activity or should it be considered as a victim of its own success and its practitioners depart to fields anew -- including new areas of activity that have been fostered or enabled by GNI?

These are difficult questions and we claim no special authority to answer them in an emphatic fashion. Yet, these are questions which, we believe, must be addressed. This short article is an attempt to foster a discussion. We commence with a brief survey of the main themes of GNI {\em circa\/} 2015. This is followed by a review of recent and ongoing developments, as well as of some new research directions that have emerged from GNI but have acquired a life of their own. 

\section{The story so far}

\subsection{Symplectic integration}

The early story of GNI is mostly the story of symplectic methods. A Hamiltonian system
\begin{equation}
  \label{Hamiltonian}
  \dot{\MM{p}}=-\frac{\partial H(\MM{p},\MM{q})}{\partial\MM{q}},\qquad \dot{\MM{q}}=\frac{\partial H(\MM{p},\MM{q})}{\partial \MM{p}},
\end{equation}
where $H:\BB{R}^{2d}\rightarrow\BB{R}$ is a {\em Hamiltonian energy,\/} plays a fundamental role in mechanics and is known to possess a long list of structural invariants, e.g.\ the conservation of the Hamiltonian energy. Yet, arguably its most important feature is the conservation of the {\em symplectic form\/} $\sum_{k=1}^d \D\MM{p}_k\wedge\D\MM{q}_k$ because symplecticity is equivalent to Hamiltonicity -- in other words, every solution of a Hamiltonian system is a symplectic flow and every symplectic flow is Hamiltonian with respect to an appropriate Hamiltonian energy \cite{hairer06gni}. 

The  solution of Hamiltonian problems using symplectic methods has a long history, beautifully reviewed in \cite{hairer03gni}, but modern efforts can be traced to the work of Feng and his collaborators at the Chinese Academy of Sciences, who have used generating-function methods to solve Hamiltonian systems \cite{feng89ccd}. And then, virtually simultaneously, \citeasnoun{lasagni88crk}, \citeasnoun{sanzserna88rks} and \citeasnoun{suris88pss} proved that certain Runge--Kutta methods, including the well-known Gauss--Legendre methods, preserve symplecticity and they presented an easy criterion for the symplecticity of Runge--Kutta methods. GNI came of age!

Subsequent research into symplectic Runge--Kutta methods had branched out into a number of directions, each with its own important ramifications outside the Hamiltonian world:
\begin{itemize}
\item {\em Backward error analysis.\/} The idea of backward error analysis (\reflectbox{BEA})  can be traced to  Wilkinson's research into linear algebra algorithms in the 1950ties. Instead of asking ``what is the numerical error for our problem", Wilkinson asked ``which nearby problem is solved {\em exactly\/} by our method?". The difference between the original and the nearby problem can tell us a great deal about the nature of the error in a numerical algorithm. 

A generalisation of \reflectbox{BEA} to the field of differential equations is fraught with difficulties. Perhaps the first successful attempt to analyse Hamiltonian ODEs in this setting was by \citeasnoun{neishtadt84sms} and it was followed by many, too numerous to list: an excellent exposition (like for many things GNI) is the monograph of \citeasnoun{hairer06gni}. The main technical tool is the B-series, an expansion of composite functions in terms of forests of rooted trees, originally pioneered by \citeasnoun{butcher63csr}. (We mention in passing that the Hopf algebra structure of this {\em Butcher group\/} has been recently exploited by mathematical physicists to understand the renormalisation group \cite{connes99lqt} -- as the authors write, ``We regard Butcher's work on the classification of numerical integration methods as an impressive example that concrete problem-oriented work can lead to far-reaching conceptual results''.) It is possible to prove that, subject to very generous conditions, the solution of a Hamiltonian problem by a symplectic method, implemented with constant step size, is exponentially near to the {\em exact\/} solution of a nearby Hamiltonian problem for an exponentially long time. This leads to considerably greater numerical precision, as well as to the conservation on average (in a strict ergodic sense) of Hamiltonian energy.

B-series fall short in a highly oscillatory and multiscale setting, encountered frequently in practical Hamiltonian systems. The alternative in the \reflectbox{BEA} context is an expansion into {\em modulated Fourier series\/} \cite{hairer00lte}.

\item {\em Composition and splitting.}

Many Hamiltonians of interest can be partitioned into a sum of kinetic and potential energy, $H(\MM{p},\MM{q})=\MM{p}^\top M\MM{p}+V(\MM{q})$.  It is often useful to take advantage of this in the design of symplectic methods. While conventional symplectic Runge--Kutta methods are implicit, hence expensive, {\em partitioned Runge--Kutta methods,\/} advancing separately in the `direction' of kinetic and potential energy, can be explicit and are in general much cheaper. While  perhaps the most important method, the St\"ormer--Verlet scheme \cite{hairer03gni}, has been known for many years, modern theory has led to an entire menagerie of composite and partitioned methods  \cite{sanzserna94nhp}. 

Splitting methods\footnote{Occasionally known in the PDE literature as {\em alternate direction methods.\/}} have been used in the numerical solution of PDEs since 1950s. Thus, given the equation $u_t=\mathcal{L}_1(u)+\mathcal{L}_2(u)$, where the $\mathcal{L}_k$s are (perhaps nonlinear) operators, the idea is to approximate the solution in the form
\begin{equation}
  \label{splitting}
  u(t+h)\approx \ee^{\alpha_1 h\mathcal{L}_1} \ee^{\beta_1 h\mathcal{L}_2} \ee^{\alpha_2 h\mathcal{L}_1} \cdots \ee^{\alpha_s h\mathcal{L}_1} \ee^{\beta_s\mathcal{L}_2}u(t),
\end{equation}
where $v(t_0+h)=:\ee^{h \mathcal{L}_1}v(t_0)$ and $w(t_0+h)=:\ee^{h \mathcal{L}_2}w(t_0)$ are, formally, the solutions of $\dot{v}=\mathcal{L}_1(v)$ and $\dot{w}=\mathcal{L}_2(w)$ respectively, with suitable boundary conditions. The underlying assumption is that the solutions of the latter two equations are either available explicitly or are easy to approximate, while the original equation is more difficult. 

A pride of place belongs to {\em palindromic compositions\/} of the form
\begin{equation}
  \label{palindromic}
  \ee^{\alpha_1 h\mathcal{L}_1} \ee^{\beta_1 h\mathcal{L}_2} \ee^{\alpha_2 h\mathcal{L}_1} \cdots \ee^{\alpha_q h\mathcal{L}_1}\ee^{\beta_q h\mathcal{L}_2}\ee^{\alpha_q h\mathcal{L}_1} \cdots \ee^{\alpha_2 h\mathcal{L}_1}  \ee^{\beta_1 h\mathcal{L}_2} \ee^{\alpha_1 h\mathcal{L}_1},
\end{equation}
invariant with respect to a reversal of the terms. They constitute a {\em time-symmetric map,\/} and this has a number of auspicious consequences. Firstly, they are always of an even order. Secondly -- and this is crucial in the GNI context -- they respect both structural invariants whose integrators are closed under composition, i.e.\ form a group (for example integrators preserving volume, symmetries, or first integrals), as well as invariants whose integrators are closed under symmetric composition, i.e.\ form a symmetric space (for example  integrators that are self-adjoint, or preserve reversing symmetries). A basic example of \R{palindromic} is the second-order {\em Strang composition\/}
\begin{displaymath}
  \ee^{\frac12 h\mathcal{L}_1} \ee^{h\mathcal{L}_2}\ee^{\frac12 h \mathcal{L}_1} =\ee^{h(\mathcal{L}_1+\mathcal{L}_2)} +\O{h^3}.
\end{displaymath}
Its order -- and, for that matter, the order of any time-symmetric method -- can be boosted by the {\em Yoshida device\/} \cite{yoshida90cho}. Let $\Phi$ be a time-symmetric approximation to $\ee^{t\mathcal{L}}$ of order $2P$, say. Then
\begin{displaymath}
  \Phi((1+\alpha)h)\Phi(-(1+2\alpha)h)\Phi((1+\alpha)h),\qquad \mbox{where}\qquad \alpha=\frac{2^{1/(2P+1)}-1}{2-2^{1/(2P+1)}}
\end{displaymath}
is also time symmetric and of order $2P+2$.  Successive applications of the Yoshida device allow to increase arbitrarily the order of the Strang composition, while retaining its structure-preserving features. This is but a single example of the huge world of splitting and composition methods, reviewed in \cite{mclachlan02sm}.

\item {\em Exponential integrators.}

Many `difficult' ODEs can be written in the form $\dot{\MM{y}}=A\MM{y}+\MM{b}(\MM{y})$ where the matrix $A$ is `larger' (in some sense) than $\MM{b}(\MM{y})$ -- for example, $A$ may be the Jacobian of an ODE (which may vary from step to step). Thus, it is to be expected that the `nastiness' of the ODE under scrutiny -- be it stiffness, Hamiltonicity or high oscillation -- is somehow `hardwired' into the matrix $A$. The exact solution of the ODE can be written in terms of the variation-of-constants formula,
\begin{equation}
  \label{VoC}
  \MM{y}(t+h)=\ee^{hA}\MM{y}(t)+\int_0^h \ee^{(h-\xi)A}\MM{b}(\MM{y}(t+\xi))\D\xi,
\end{equation}
except that, of course, the right-hand side includes the unknown function $\MM{y}$. Given the availability of very effective methods to compute the matrix exponential, we can exploit this to construct {\em exponential integrators,\/} explicit methods that often exhibit favourable stability and structure-preservation features. The simplest example, the {\em exponential Euler\/} method, freezes $\MM{y}$ within the integral in \R{VoC} at its known value at $t$, the outcome being the first-order method
\begin{displaymath}
  \MM{y}_{n+1}=\ee^{hA}\MM{y}_n+A^{-1}(\ee^{hA}-I)\MM{b}(\MM{y}_n).
\end{displaymath}
The order can be boosted by observing that (in a loose sense which can be made much more precise) the integral above is discretised by the Euler method, which is a one-stage explicit Runge--Kutta scheme, discretising it instead by multistage schemes of this kind leads to higher-order methods \cite{hochbruck10ei}. 

Many Hamiltonian systems of interest can be formulated as second-order systems of the form $\ddot{\MM{y}}+\Omega^2\MM{y}=\MM{g}(\MM{y})$. Such systems feature prominently in the case of highly oscillatory mechanical systems, where $\Omega$ is positive definite and has some large eigenvalues. Variation of constants \R{VoC} now reads
\begin{Eqnarray*}
    \left[\!
  \begin{array}{c}
     \MM{y}(t+h)\\
     \dot{\MM{y}}(t+h)
  \end{array}
  \!\right]&=&  \left[
  \begin{array}{cc}
     \cos(h\Omega) & \Omega^{-1}\sin(h\Omega)\\
     -\Omega\sin(h\Omega) & \cos(h\Omega)
  \end{array}
  \right]  \left[\!
  \begin{array}{c}
     \MM{y}(t)\\
     \dot{\MM{y}}(t)
  \end{array}
  \!\right]\\
  &&\mbox{}+\int_t^{t+h}   \left[
  \begin{array}{cc}
     \cos((h-\xi)\Omega) & \Omega^{-1}\sin((h-\xi)\Omega)\\
     -\Omega\sin((h-\xi)\Omega) & \cos((h-\xi)\Omega)
  \end{array}
  \right] \! \left[
  \begin{array}{c}
     \MM{0}\\
     \MM{g}(\MM{y}(t+\xi))
  \end{array}
  \right]\!\D\xi
\end{Eqnarray*}
and we can use either standard exponential integrators or exponential integrators designed directly for second-order systems and using  Runge--Kutta--Nystr\"om methods on the nonlinear part \cite{wu13spa}. 

An important family of exponential integrators for second-order systems are {\em Gautschi-type methods\/}
\begin{equation}
  \label{Gautschi}
  \MM{y}_{n+1}-2\MM{y}_n+\MM{y}_{n-1}=h^2\Psi(h\Omega) (\MM{g}_n-\Omega^2\MM{y}_n),
\end{equation}
which are of second order. Here $\Psi(x)=2(1-\cos x)/x$ while, in Gautschi's original method, $\MM{g}_n=\MM{g}(\MM{y}_n)$ \cite{hochbruck10ei}. Unfortunately, this choice results in resonances and a better one is $\MM{g}_n=\MM{g}(\Phi(h\Omega)\MM{y}_n)$, where the {\em filter\/} $\Phi$ eliminates resonances: $\Phi(0)=I$ and $\Phi(k\pi)=0$ for $k\in\BB{N}$. We refer to \cite{hochbruck10ei} for further discussion of such methods in the context of symplectic integration. 

\item {\em Variational integrators.}  {\em Lagrangian formulation\/} recasts a large number of differential equations as minima of nonlinear functionals. Thus, for example, instead of the Hamiltonian problem $M\ddot{\MM{q}}+\MM{\nabla} V(\MM{q})=\MM{0}$, where the matrix $M$ is positive definite, we may consider the equivalent variational formulation of minimizing the positive-definite nonlinear functional $L(\MM{q},\dot{\MM{q}})=\frac12 \dot{\MM{q}}^\top M\dot{\MM{q}}-V(\MM{q})$. With greater generality, Hamiltonian and Lagrangian formulations are connected via the familiar Euler--Lagrange equations and, given the functional $L$, the corresponding second-order  system is
\begin{displaymath}
  \frac{\partial L(\MM{q},\dot{\MM{q}})}{\partial\MM{q}}-\frac{\D}{\D t} \left[\frac{\partial L(\MM{q},\dot{\MM{q}})}{\partial \dot{\MM{q}}}\right]=\MM{0}.
\end{displaymath}
The rationale of variational integrators parallels that of the {\em Ritz method\/} in the theory of finite elements. We first reformulate the Hamiltonian problem as a Lagrangian one, project it to a finite-dimensional space, solve it there and transform back. The original symplectic structure is replaced by a finite-dimensional symplectic structure, hence the approach is by design symplectic \cite{marsden01dmv}. 
\end{itemize}

\subsection{Lie-group methods}

Let $\mathcal{G}$ be a Lie group and $\mathcal{M}$ a differentiable manifold. We say that $\Lambda:\mathcal{G}\times\mathcal{M}\rightarrow\mathcal{M}$ is a {\em group action\/} if\\[4pt]
a.~$\Lambda(\iota,y)=y$ for all $y\in\mathcal{M}$ (where $\iota$ is the identity of $\mathcal{G}$) and \\[2pt]
b.~$\Lambda(p,\Lambda(q,y))=\Lambda(p\cdot q,y)$ for all $p,q\in\mathcal{G}$ and $y\in\mathcal{M}$.\\[4pt]
If, in addition, for every $x,y\in\mathcal{M}$ there exists $p\in\mathcal{G}$ such that $y=\Lambda(p,x)$, the action is said to be transitive and $\mathcal{M}$ is a {\em homogeneous space,\/} acted upon by $\mathcal{G}$. 

Every Lie group acts upon itself, while the orthogonal group $\CC{O}(n)$  acts on the $(n-1)$-sphere  by multiplication, $\Lambda(p,y)=py$. The orthogonal group also acts on the {\em isospectral manifold\/} of all symmetric matrices similar to a specific symmetric matrix by similarity, $\Lambda(p,y)=pyp^\top$. Given $1\leq m\leq n$, the {\em Grassmann manifold\/} $\BB{G}(n,m)$ of all $m$-dimensional subspaces of $\BB{R}^n$ is a homogeneous space acted upon by $\CC{SO}(m)\times\CC{SO}(n-m)$, where $\CC{SO}(m)$ is the special orthogonal group -- more precisely, $\BB{G}(n,m)=\CC{SO}(n)/(\CC{SO}(m)\times\CC{SO}(n-m))$. 

Faced with a differential equation evolving in a homogeneous space, we can identify its flow with a group action: Given an initial condition $y_0\in\mathcal{M}$, instead of asking ``what is the value of $y$ at time $t>0$'' we might pose the equivalent question ``what is the group action that takes the solution from $y_0$ to $y(t)$?''. This is often a considerably more helpful formulation because a group action can be further related to an {\em algebra action.\/} Let $\GG{g}$ be the Lie algebra corresponding to the  matrix group $\mathcal{G}$, i.e.\ the tangent space at $\iota\in\mathcal{G}$, and denote by $\GG{X}(\mathcal{M})$ the set of all Lipschitz vector fields over $\mathcal{M}$. Let $\lambda:\GG{g}\rightarrow\GG{X}(\mathcal{M})$ and $a:\BB{R}_+\times\mathcal{M}\rightarrow\GG{g}$ be both Lipschitz. In particular, we might consider 
\begin{displaymath}
  \lambda(a,y)=\frac{\D}{\D s} \Lambda(\rho(s,y),y)\,\rule[-6pt]{0.75pt}{18pt}_{\,s=0},
\end{displaymath}
where $\Lambda$ is a group action and $\rho:\BB{R}_+\rightarrow\mathcal{G}$, $\rho(s,y(s))=\iota+a(s,y(s))s+\O{s^2}$ for small $|s|$. The equation $\dot{y}=\lambda(a(t,y),y)$, $y(0)=y_0\in\mathcal{M}$ represents {\em algebra action\/} and its solution evolves in $\mathcal{M}$. Moreover, 
\begin{equation}
  \label{alg_action}
  y(t)=\Lambda(v(t),y_0)\qquad \mbox{where}\qquad \dot{v}=a(t,\Lambda(v,y_0))v,\quad v(0)=\iota\in\mathcal{G}
\end{equation}
is a {\em Lie-group equation.\/} Instead of solving the original ODE on $\mathcal{M}$, it is possible to solve \R{alg_action} and use the group action $\Lambda$ to advance the solution to the next step: this is the organising principle of most {\em Lie-group methods\/}  \cite{iserles00lgm}. It works because a Lie-group equation can be solved in the underlying Lie algebra, which is a {\em linear space.\/} Consider an ODE\footnote{Or, for that matter, a PDE, except that formalities are somewhat more complicated.} $\dot{y}=f(y)$, $y(0)\in\mathcal{M}$, such that $f:\mathcal{M}\rightarrow\GG{X}$ -- the solution $y(t)$ evolves on the manifold. While conventional numerical methods are highly unlikely to stay in $\mathcal{M}$, this is not the case for Lie-group methods. We can travel safely  between $\mathcal{M}$ and $\mathcal{G}$ using a group action. The traffic between $\mathcal{G}$ and $\GG{g}$ is slightly more complicated and we need to define a {\em trivialisation,\/} i.e.\ an invertible map taking smoothly a neighbourhood of $0\in\GG{g}$ to a neighbourhood of $\iota\in\mathcal{G}$ and taking zero to identity. The most ubiquitous example of trivialisation is the exponential map, which represents the solution of \R{alg_action} as $v(t)=\ee^{\omega(t)}$, where $\omega$ is the solution of the {\em dexpinv equation\/}
\begin{equation}
  \label{dexpinv}
  \dot{\omega}=\sum_{m=0}^\infty \frac{\CC{B}_m}{m!} \CC{ad}^m_{a(t,\ee^\omega)}\omega,\qquad \omega(0)=0\in\GG{g}
\end{equation}
\cite{iserles00lgm}. Here the $\CC{B}_m$s are Bernoulli numbers, while $\CC{ad}^m_b$ is the {\em adjoint operator\/} in $\GG{g}$,
\begin{displaymath}
  \CC{ad}_b^0 c=c,\qquad \CC{ad}_b^m c=[b,\CC{ad}_b^{m-1}c],\quad m\in\BB{N},\qquad b,c\in\GG{g}.
\end{displaymath}
Because $\GG{g}$ is closed under linear operations and commutation, solving \R{dexpinv} while respecting Lie-algebraic structure is straightforward. Mapping back, first to $\mathcal{G}$ and finally to $\mathcal{M}$, we keep the numerical solution of $\dot{y}=f(t)$ on the manifold. 

Particularly effective is the use of explicit Runge--Kutta methods for \R{dexpinv}, the so-called Runge--Kutta--Munthe-Kaas (RKMK) methods \cite{munthekass98rkm}. To help us distinguish between conventional Runge--Kutta methods and RKMK, consider the three-stage, third-order method with the Butcher tableau\footnote{For traditional concepts such as Butcher tableaux, Runge-Kutta methods and B-series, the reader is referred to \cite{hairer93sod}.}
\begin{equation}
  \label{RK3}
  \begin{array}{c|ccc}
     0 & \\
     \frac12 & \frac12\\[2pt]
     1 & -1 & 2\\\hline
     & \frac16 & \frac23 & \frac16\rule{0pt}{13pt}
  \end{array}.
\end{equation}
Applied to the  ODE $\dot{y}=f(t,y)$, $y(t_n)=y_n\in\mathcal{M}$, evolving on the manifold $\mathcal{M}\subset\BB{R}^d$, it becomes
\begin{Eqnarray*}
  &&k_1=f(t_n,y_n),\\
  && k_2=f(t_{n+\frac12},y_n+\Frac12 hk_1),\\
  && k_3=f(t_{n+1},y_n-hk_1+2hk_2),\\
  &&\Delta=h(\Frac16 k_1+\Frac23 k_2+\Frac16 k_3),\\[3pt]
  y_{n+1}&=&y_n+\Delta.
\end{Eqnarray*}
Since we operate in $\BB{R}^d$, there is absolutely no reason for $y_{n+1}$ to live in $\mathcal{M}$.  However, once we implement \R{RK3} at an algebra level (truncating first the dexpinv equation \R{dexpinv}),
\begin{Eqnarray*}
  &&k_1=a(t_n,\iota),\\
  &&k_2=a(t_{n+\frac12},\ee^{hk_1/2}),\\
  &&k_3=a(t_{n+1},\ee^{-hk_1+2hk_2}),\\
  &&\Delta=h(\Frac16 k_1+\Frac23 k_2+\Frac16 k_3),\\[3pt]
  \omega_{n+1}&=&\Delta+\Frac16 h[\Delta,k_1]\\
  y_{n+1}&=&\Lambda(\ee^{\omega_{n+1}},y_n),
\end{Eqnarray*}
the solution is guaranteed to stay in $\mathcal{M}$. 

An important special case of a Lie-group equation is the linear ODE $\dot{v}=a(t)v$, where $a:\BB{R}_+\rightarrow\GG{g}$. Although RKMK works perfectly well in a linear case, special methods do even better. Perhaps the most important is the {\em Magnus expansion\/} \cite{magnus54esd}, $v(t)=\ee^{\omega(t)}v(0)$, where
\begin{Eqnarray}
  \nonumber
  \omega(t)&=& \int_0^t a(\xi)\D \xi -\frac12 \int_0^t\!\int_0^{\xi_1} [a(\xi_2),a(\xi_1)]\D\xi_2\D\xi_1 \\
  \label{Magnus}
  &&\mbox{}+\frac14 \int_0^t \! \int_0^{\xi_1} \! \!\int_0^{\xi_2} [[a(\xi_3),a(\xi_2)],a(\xi_1)]\D\xi_3\D\xi_2\D\xi_1\\
  \nonumber
  &&+\frac{1}{12} \int_0^t\!\int_0^{\xi_1}\!\!\int_0^{\xi_2}[a(\xi_3),[a(\xi_2),a(\xi_1)]]\D\xi_3\D\xi_2\D\xi_1+\cdots.
\end{Eqnarray}
We refer to \cite{iserles99sld,iserles00lgm,blanes09mes} for explicit means to derive expansion terms, efficient computation of multivariate integrals that arise in this context and many other implementation details. Magnus expansions are important in a number of settings when preservation of structure is not an issue, not least in the solution of linear stochastic ODEs \cite{lord08esi}.
 
There are alternative means to expand the solution of \R{dexpinv} in a linear case, not least the {\em Fer expansion,\/} that has found recently an important application in the computation of Sturm--Liouville spectra \cite{ramos15nss}. 

Another approach to Lie-group equations uses {\em canonical coordinates of the second kind\/} \cite{owren01imb}.

\subsection{Conservation of volume}

An ODE $\dot{\MM{x}}=\MM{f}(\MM{x})$ is divergence-free if $\MM{\nabla} \cdot \MM{f}(\MM{x})=0$. The flows of divergence-free ODEs are volume-preserving (VP). Volume  is important to preserve, as it leads to KAM-tori, incompressibility, and, most importantly, is a crucial ingredient for ergodicity. Unlike symplecticity, however, phase space volume can generically {\it not} be preserved by Runge--Kutta methods, or even by their generalisations, B-series methods. This was proved independently in \cite{chartier07pfi} and in \cite{iserles07bmc}. Since B-series methods cannot preserve volume, we need to look to other methods.

There are essentially two known numerical integration methods that preserve phase space volume. The first volume-preserving method is based on splitting  \cite{feng95vpa}. As an example, consider a 3D volume preserving vector field:
\begin{Eqnarray}
\dot{x} &=& u(x,y,z) \nonumber \\
\label{3D}
\dot{y} &=& v(x,y,z)  \\ 
\dot{z} &=& w(x,y,z)  \nonumber 
\end{Eqnarray}
with
\begin{displaymath} 
u_x + v_y + w_z = 0. 
\end{displaymath}
We split this 3D VP vector field into two 2D VP vector fields as follows
\begin{equation}
  \label{VP1}
  \begin{array}{lcl}
    \displaystyle \dot{x} = u(x,y,z), &\qquad\quad& \displaystyle \dot{x} = 0,\\[6pt]
    \displaystyle \dot{y} = -\int\! u_x(x,y,z)\D y, && \displaystyle \dot{y} = v(x,y,z) + \int\! u_x(x,y,z) \D y,\\[12pt]
    \displaystyle \dot{z} = 0; && \displaystyle \dot{z} = w(x,y,z).
  \end{array}\hspace*{20pt}
\end{equation}
The vector field on the left is divergence-free by construction, and since both vector fields add up to (2.1), it follows that the vector field on the right is also volume-preserving.

Having split the original vector field into 2D VP vector fields, we need to find VP integrators for each of these 2D VP vector fields. But that is easy, since 2D VP vector fields are essentially equivalent to 2D Hamiltonian vector fields (with the extra dimension `frozen'), and all symplectic methods (e.g. symplectic Runge--Kutta methods) are volume-preserving for Hamiltonian vector fields.

The above splitting method is easily generalised to $n$ dimensions, where one splits into $n-1$ 2D VP vector fields, and integrates each  using a symplectic Runge--Kutta method.

An alternative VP integration method was discovered independently by Shang and by Quispel \cite{shang94gfv,quispel95vpi}. We again illustrate this method in 3D.

We will look for an integrator of the form
\begin{Eqnarray}
x_1 &=& g_1(x_1',x_2,x_3) \nonumber \\
\label{VPint}
x_2' &=& g_2(x_1',x_2,x_3)  \\ 
x_3' &=& g_1(x_1',x_2',x_3)  \nonumber 
\end{Eqnarray}
where (here and below) $x_i= x_i(nh)$, and $x_i'=x_i((n+1)h)$. The reason the form \R{VPint} is convenient, is because any such map is VP iff
\begin{equation}\label{VPintcon}
\frac{\partial x_1}{\partial x_1'} = \frac{\partial x_2'}{\partial x_2}\frac{\partial x_3'}{\partial x_3}.
\end{equation}
To see how to construct a VP integrator of the form \R{VPint}, consider as an example the ODE
\begin{Eqnarray}
\dot{x}_1 &=& x_2 + x_1^2 + x_3^3 \nonumber \\
\label{3D2}
\dot{x}_2 &=& x_3 + x_1x_2 + x_1^4  \\ 
\dot{x}_3 &=& x_1 - 3x_1x_3 + x_2^5  \nonumber 
\end{Eqnarray}
It is easy to check that it is divergence-free.

\bigskip \noindent Now consistency requires that any integrator for \R{3D2} should satisfy
\begin{Eqnarray}
x_1' &=& x_1 + h( x_2 + x_1^2 + x_3^3) + \O{h^2} \nonumber \\
\label{consis1}
x_2' &=& x_2 + h(x_3 + x_1x_2 + x_1^4) + \O{h^2}  \\ 
x_3' &=& x_3 + h(x_1 - 3x_1x_3 + x_2^5) + \O{h^2}  \nonumber 
\end{Eqnarray}
and therefore
\bigskip \noindent 
\begin{Eqnarray}
x_1 &=& x_1' - h( x_2 + x_1'^2 + x_3^3) + \O{h^2}  \\
\label{consis2}
x_2' &=& x_2 + h(x_3 + x_1'x_2 + x_1'^4) + \O{h^2}  \\ 
x_3' &=& x_3 + h(x_1' - 3x_1'x_3 + x_2'^5) + \O{h^2} 
\end{Eqnarray}
Since we are free to choose any consistent $g_2$ and $g_3$ in \R{VPint}, provided $g_1$ satisfies \R{VPintcon}, we choose the terms designated by $\O{h^2}$ in (2.15) and (2.16) to be identically zero. Equation \R{VPintcon} then yields
\begin{equation}\label{VPegcon}
\frac{\partial x_1}{\partial x_1'} = (1+hx_1')(1-3hx_1').
\end{equation}
This can easily be integrated to give
\begin{equation}\label{VPintconsol}
x_1 = x_1' - hx_1'^2 - h^2x_1'^3 + k(x_2,x_3;h).
\end{equation}
where the function $k$ denotes an integration constant that we can choose appropriately. The simplest VP integrator satisfying both (2.14) and \R{VPintconsol} is therefore:
\begin{Eqnarray}
x_1 &=& x_1' - h( x_2 + x_1'^2 + x_3^3) -h^2x_1'^3 \nonumber  \\
\label{egint}
x_2' &=& x_2 + h(x_3 + x_1'x_2 + x_1'^4)   \\ 
x_3' &=& x_3 + h(x_1' - 3x_1'x_3 + x_2'^5) \nonumber
\end{Eqnarray}
A nice aspect of the integrator \R{egint}  (and \R{VPint}) is that it is essentially only implicit in one variable. Once $x_1'$ is computed from the first (implicit) equation, the other two equations are essentially explicit.

\bigskip \noindent Of course the method just described also generalises to any divergence-free ODE in any dimension.

\subsection{Preserving energy and other first integrals}

As mentioned, Hamiltonian systems exhibit two important geometric properties simultaneously, they conserve both the symplectic form and the energy. A famous no-go theorem by \citeasnoun{ge88lph} has shown that it is generically impossible to construct a geometric integrator that preserves both properties at once. One therefore must choose which one of these two to preserve in any given application. Particularly in low dimensions and if the energy surface is compact, there are often advantages in preserving the energy.

An energy-preserving B-series method was discovered in \cite{quispel08nce} cf.\ also \cite{mclachlan99gid}.

For any ODE $\dot{\MM{x}} = \MM{f}(\MM{x})$, this so-called average vector field method is given by
\begin{equation}\label{defavf}
\frac{\MM{x}'-\MM{x}}{h} = \int_{0}^{1}\MM{f}(\xi \MM{x}' + (1-\xi)\MM{x})\D\xi.
\end{equation}
If the vector field $\MM{f}$ is Hamiltonian, i.e. if there exists a Hamiltonian function $H(\MM{x})$ and a constant skew-symmetric matrix  $S$ such that $\MM{f}(\MM{x}) = S\nabla H(\MM{x})$, then it follows from \R{defavf} that energy is preserved, i.e. $H(\MM{x}')=H(\MM{x})$.

While the B-series method \R{defavf} is energy-preserving for any Hamiltonian $H$,  it can be shown that no Runge--Kutta method is energy-preserving for all $H$. For a given {\it polynomial} $H$ however, Runge--Kutta methods preserving that $H$ do exist \cite{iavtrig09hos}. This can be seen as follows.

Note that the integral in \R{defavf} is one-dimensional. This means that e.g.\ for cubic vector fields (and hence for quartic Hamiltonians) an equivalent method is obtained by replacing the integral in \R{defavf} using Simpson's rule:
\begin{equation}\label{simp}
 \int_{0}^{1} g(\xi)\D\xi \approx \frac{1}{6}\left[g(0) + 4g(\Frac{1}{2}) + g(1)\right]\!.
\end{equation}
yielding the Runge--Kutta method
\begin{equation}\label{RKsimp}
 \frac{\MM{x}'-\MM{x}}{h} =  \frac{1}{6}\left[\MM{f}(\MM{x}) + 4\MM{f}\left(\frac{\MM{x}+\MM{x}'}{2}\right) + \MM{f}(\MM{x}')\right]\!,
\end{equation}
preserving all quartic Hamiltonians.

We note that \R{defavf} has second order accuracy. Higher order generalisations have been given in \cite{hairer10epv}. We note that the average vector field method has also been applied to a slew of semi-discretised PDEs in \cite{celledoni12per}.

While energy is one of the most important constants of the motion in applications, many other types of first integrals do occur. We note here that all B-series methods preserve all linear first integrals, and that all symplectic B-series methods preserve all quadratic first integrals. So, for example, the implicit midpoint rule 
\begin{displaymath}
  \frac{\MM{x}'-\MM{x}}{h} = \MM{f}\!\left(\frac{\MM{x}+\MM{x}'}{2}  \right)
\end{displaymath}
(which is symplectic) preserves all linear and quadratic first integrals. There are however many cases not covered by any of the above.

How does one preserve a cubic first integral that is not energy?  And what about Hamiltonian systems whose symplectic structure is not constant? It turns out that generically, any ODE $\dot{\MM{x}} = \MM{f}(\MM{x})$ that preserves an integral $I(\MM{x})$, can be written in the form
\begin{equation}\label{gendeform}
 \dot{\MM{x}} = S(\MM{x})\MM{\nabla} I(\MM{x}),
\end{equation}
where $S(\MM{x})$ is a skew-symmetric matrix\footnote{Note that in general $S(\MM{x})$ need not satisfy the so-called Jacobi identity.}.

An integral-preserving discretisation of \R{gendeform} is given by
\begin{equation}
  \label{gendeintgrtr}
 \frac{\MM{x}'-\MM{x}}{h} =  \bar{S}(\MM{x},\MM{x}') \bar{\nabla}I(\MM{x},\MM{x}'),
\end{equation}
where $\bar{S}(\MM{x},\MM{x}')$ is any consistent approximation to $S(\MM{x})$ (e.g. $\bar{S}(\MM{x},\MM{x}')=S(\MM{x})$), and the {\em discrete gradient\/} $ \bar{\MM{\nabla}}I$ is defined by
\begin{equation}
  \label{dgdefn1}
  (\MM{x}'-\MM{x}) \cdot \bar{\MM{\nabla}}I(\MM{x},\MM{x}') = I(\MM{x}') - I(\MM{x})
\end{equation}
and
\begin{equation}
  \label{dgdefn2}
  \lim_{\Mm{x}' \rightarrow \Mm{x}}  \bar{\MM{\nabla}}I(\MM{x},\MM{x}') = \MM{\nabla} I(\MM{x}).
\end{equation}
There are many different discrete gradients that satisfy \R{dgdefn1} and \R{dgdefn2}. A particularly simple one is given by the Itoh--Abe discrete gradient, which for example in 3D reads
\begin{equation}\label{ItohAbe}
\bar{\nabla}I(\MM{x},\MM{x}') = \left[
\begin{array}{c}
\displaystyle \frac{I(x_1',x_2,x_3) - I(x_1,x_2,x_3)}{x_1'-x_1} \\[10pt]
\displaystyle \frac{I(x_1',x_2',x_3) - I(x_1',x_2,x_3)}{x_2'-x_2} \\[10pt]
\displaystyle   \frac{I(x_1',x_2',x_3') - I(x_1',x_2',x_3)}{x_3'-x_3} 
  \end{array}
\right]\!.
\end{equation}
Other examples of discrete gradients, as well as constructions of the skew-symmetric matrix $S(\MM{x})$ for a given vector field $\MM{f}$ and integral $I$ may be found in \cite{mclachlan99gid}.

We note that the discrete gradient method can also be used for systems with any number of integrals. For example an ODE $\dot{\MM{x}}=\MM{f}(\MM{x})$ possessing two integrals $I(\MM{x})$ and $J(\MM{x})$ can be written
\begin{equation}\label{ode2ints}
\dot{x}_i = S_{ijk}(\MM{x}) \frac{\partial I(\MM{x})}{\partial x_j} \frac{\partial J(\MM{x})}{\partial x_k},
\end{equation}
where the summation convention is assumed over repeated indices and $S(\MM{x})$ is a completely antisymmetric tensor. A discretisation of \R{ode2ints} which preserves both $I$ and $J$ is given by
\begin{equation}\label{disc2ints}
\frac{x_i' - x_i}{h} = \bar{S}_{ijk}(\MM{x},\MM{x}') \bar{\nabla}I(\MM{x},\MM{x}')  \,\rule[-4pt]{0.5pt}{16pt}_{\,j} \bar{\nabla}J(\MM{x},\MM{x}') \,\rule[-4pt]{0.5pt}{16pt}_{\,k}
\end{equation}
with $\bar{S}$ any completely skew approximation of $S$ and $\bar{\nabla}I$ and $\bar{\nabla}J$ discrete gradients as defined above.

\setcounter{equation}{0}
\setcounter{figure}{0}
\section{Five recent stories of GNI}

The purpose of this section is not to present a totality of recent research into GNI, a subject that would have called for a substantially longer paper. Instead, we wish to highlight a small number of developments with which the authors are familiar and which provide a flavour of the very wide range of issues on the current GNI agenda.

\subsection{Highly oscillatory Hamiltonian systems} 

High oscillation occurs in many Hamiltonian systems. Sometimes, e.g.\ in the integration of equations of celestial mechanics, the source of the problem is that we wish to compute the solution across a very large number of periods and the oscillation is an artefact of the time scale in which the solution has physical relevance. In other cases oscillation is implicit in the multiscale structure of the underlying problem. A case in point are the (modified) {\em Fermi--Pasta--Ulam (FPU) equations,\/} describing a mechanical system consisting of alternating stiff harmonic and soft nonlinear springs. The soft springs impart fast oscillation, while the hard springs generate slow transfer of energy across the system: good numerical integration must capture both!

A good point to start (which includes modified FPU as a special case) is the second-order ODE
\begin{equation}
  \label{HiOscODE}
  \ddot{\MM{q}}+\Omega^2\MM{q}=\MM{g}(\MM{q}),\qquad t\geq0,\qquad \MM{q}(0)=\MM{u}_0,\quad \dot{\MM{q}}(0)=\MM{v}_0,
\end{equation}
where $\MM{g}(\MM{q})=-\MM{\nabla}U(\MM{q})$ and
\begin{displaymath}
  \Omega=  \left[
  \begin{array}{cc}
     O & O\\
     O & \omega I
  \end{array}
  \right]\!,\quad \omega\gg1,\qquad \MM{q}=  \left[
  \begin{array}{c}
     \MM{q}_0\\\MM{q}_1
  \end{array}
  \right]\!,\qquad \MM{q}_0\in\BB{R}^{n_0},\;\;\MM{q}_1\in\BB{R}^{n_1}.
\end{displaymath}
An important aspect of systems of the form \R{HiOscODE} is that the exact solution, in addition to preserving the total Hamiltonian energy
\begin{equation}
  \label{HamEn}
  H(\MM{p},\MM{q})=\frac12 (\|\MM{p}_1\|^2+\omega^2 \|\MM{q}_1\|^2)+\frac12 \|\MM{p}_0\|^2  +U(\MM{q}_0,\MM{q}_1),
\end{equation}
where $\dot{\MM{q}}=\MM{p}$, also preserves the {\em oscillatory energy\/}
\begin{equation}
  \label{OscEn}
  I(\MM{p},\MM{q})=\frac12 \|\MM{p}_1\|^2+\frac{\omega^2}{2} \|\MM{q}_1\|^2
\end{equation}
for intervals of length $\O{\omega^N}$ for any $N\geq1$. This has been proved using the {\em modulated Fourier expansions\/}
\begin{displaymath}
  \MM{q}(t)=\sum_{m=-\infty}^\infty \ee^{\ii m\omega t} \MM{z}_m(t).
\end{displaymath}

The solution of \R{HiOscODE} exhibits oscillations at frequency $\O{\omega}$ and this inhibits the efficiency of many symplectic methods, requiring step size of $\O{\omega^{-1}}$, a situation akin to stiffness in more conventional ODEs. However, by their very structure, exponential integrators (and in particular Gautschi-type methods \R{Gautschi}) are particularly effective in integrating the linear part, which gives rise to high oscillation. The problem with Gautschi-type methods, though, might be the occurrence of resonances and we need to be careful to avoid them, both in the choice of the right filter (cf.\ the discussion in Subsection~2.1) and step size $h$. 

Of course, one would like geometric numerical integrators applied to \R{HiOscODE} to exhibit favourable preservation properties with respect to both total energy \R{HamEn} and oscillatory energy \R{OscEn}. Applying modulated Fourier expansions to trigonometric and modified trigonometric integrators, this is indeed the case provided that the step size obeys the {\em non-resonance condition\/} with respect to the frequency $\omega$,
\begin{displaymath}
  |\sin(\Frac12 mh\omega)|\geq c h^{1/2},\qquad m=1,\ldots,N,\quad N\geq2,
\end{displaymath}
cf.\ \citeasnoun{hairer09olt}.

All this has been generalised to systems with multiple frequencies, with the Hamiltonian function
\begin{displaymath}
  H(\MM{p},\MM{q})=\overbrace{\frac12 \sum_{j=1}^s \left(\|\MM{p}_j\|^2+\omega_j^2\|\MM{q}_j\|^2\right)}^{\CC{oscillatory}}+\overbrace{\frac12\|\MM{p}_0\|^2+U(\MM{q})}^{\CC{slow}},
\end{displaymath}
where 
\begin{displaymath}
  \MM{p}=  \left[
  \begin{array}{c}
     \MM{p}_0\\\MM{p}_1\\\vdots\\\MM{p}_s
  \end{array}
  \right]\!,\quad \MM{q}=\left[
  \begin{array}{c}
     \MM{q}_0\\\MM{q}_1\\\vdots\\\MM{q}_s
  \end{array}
  \right]\!,\qquad 0<\min_{j=1,\ldots,s}\omega_j,\quad 1\ll \max_{j=1,\ldots,s}\omega_j
\end{displaymath}
for both the exact solution \cite{gauckler13eso} and for discretisations obtained using trigonometric and modified trigonometric integrators \cite{cohen15lta}.

Further achievements and open problem in the challenging area of marrying symplectic integration and high oscillation are beautifully described in \cite{hairer14cgi}.

\subsection{Kahan's `unconventional' method}
A novel discretisation method for quadratic ODEs was introduced and studied in \cite{kahan93unm}. This new method discretised the vector field
\begin{equation}\label{kahan1}
\dot{x}_i = \sum_{j,k}^{}a_{ijk}x_jx_k + \sum_{j}^{}b_{ij}x_j + c_i
\end{equation}
as follows,
\begin{equation}\label{kahan2}
\frac{x_i'-x_i}{h} = \sum_{j,k}^{}a_{ijk} \left(\frac{x_jx_k' + x_j'x_k}{2}\right) + \sum_{j}^{}b_{ij} \left( \frac{x_j + x_j'}{2} \right) + c_i.
\end{equation}
Kahan called the method \R{kahan2} `unconventional', because it treats the quadratic terms different from the linear terms. He also noted some nice features of \R{kahan2}, e.g. that it often seemed to be able to integrate through singularities.

\bigskip \noindent \textbf{Properties of Kahan's method:}
\begin{enumerate}
\item {\it Kahan's method is (the reduction of) a Runge--Kutta method.} 

\citeasnoun{celledoni13gpk} showed that \R{kahan2} is the reduction to quadratic vector fields of the Runge--Kutta method
\begin{equation}\label{kahan3}
\frac{\MM{x}'-\MM{x}}{h} = 2 \MM{f}\left(\frac{\MM{x} + \MM{x}'}{2}\right) - \frac{1}{2} \MM{f}(\MM{x})  - \frac{1}{2} \MM{f}(\MM{x}')
\end{equation}
This explains {\em inter alia\/} why Kahan's method preserves all linear first integrals.

\item {\it Kahan's method preserves a modified energy and measure.}

For any Hamiltonian vector field of the form
\begin{equation}\label{hamode}
\dot{\MM{x}} = \MM{f}(x) = S\MM{\nabla} H(\MM{x}),
\end{equation}
with cubic Hamiltonian $H(\MM{x})$ and constant symplectic (or Poisson) structure $S$, Kahan's method preserves a modified energy as well as a modified measure exactly \cite{celledoni13gpk}.

The modified volume is
\begin{equation}\label{modvol}
\frac{\D x_1 \wedge \dots \wedge \D x_n}{\det \!\left( I - \frac{1}{2}hf'(\MM{x}) \right)},
\end{equation}
while the modified energy is
\begin{equation}\label{modenergy}
\tilde{H}(\MM{x}) := H(\MM{x}) + \frac{1}{3}h \MM{\nabla} H(\MM{x})^\top \!\left(I - \frac{1}{2}hf'(\MM{x})  \right)^{-1} \!\MM{f}(\MM{x}).
\end{equation}

\item {\it Kahan's method preserves the integrability of many integrable systems of quadratic ODEs.}

Beginning with the work of Hirota and Kimura, subsequently extended by Suris and collaborators \cite{petrera11ihk},  and by  Quispel and collaborators \cite{celledoni13gpk,celledoni14ipk,vanderkamp14iss}, it was shown that Kahan's method preserves the complete integrability of a surprisingly large number of quadratic ODEs.

\end{enumerate}

\bigskip \noindent Here we list some 2D vector fields whose integrability is preserved  by Kahan's method:

\begin{itemize}
\item Quadratic Hamiltonian systems in 2D:

\noindent The 9-parameter family
\begin{equation}\label{9paramfam}
\left[ \begin{array}{c}
\dot{x} \\ \dot{y} \end{array} \right] =
\left[ \begin{array}{c}
bx^2 + 2cxy +dy^2 +fx + gy + i \\ -ax^2 - 2bxy - cy^2 - ex -fy -h \end{array} \right]\!;
\end{equation}

\item Suslov systems in 2D:

\noindent The 9-parameter family
\begin{equation}\label{suslov}
\left[ \begin{array}{c}
\dot{x} \\ \dot{y} \end{array} \right] = l(x,y)
\left[ \begin{array}{cc} 0 & 1 \\ -1 & 0 
\end{array} \right] \nabla H(x,y),
\end{equation}
where $l(x,y) = ax+by+c$; $H(x,y) = dx^2 + exy +fy^2 + gx + hy + i$;

\item Reduced Nahm equations in 2D:

\noindent Octahedral symmetry:
\begin{equation}\label{nahm1}
\left[\begin{array}{c}
\dot{x} \\ \dot{y} \end{array} \right] =
\left[ \begin{array}{c}
2x^2 - 12y^2 \\ -6x^2  - 4y^2  \end{array} \right]\!;
\end{equation}
Icosahedral symmetry:
\begin{equation}\label{nahm2}
\left[ \begin{array}{c}
\dot{x} \\ \dot{y} \end{array} \right] =
\left[\begin{array}{c}
2x^2 - y^2 \\ -10xy  + y^2  \end{array} \right]\!.
\end{equation}
\end{itemize}
The modified energy and measure for the Kahan discretisations of these 2D systems, as well as of many other (higher-dimensional) integrable quadratic vector fields are given in \cite{petrera11ihk,celledoni13gpk,celledoni14ipk}.

Generalisations to higher degree polynomial equations using polarisation are presented in \cite{celledoni15dpv}.

\subsection{Applications to celestial mechanics} 
GNI methods particularly come into their own when the integration time is large compared to typical periods of the system. Thus long-term integrations of e.g. solar-type systems and of particle accelerators typically need symplectic methods. In this subsection we focus on the former\footnote{A very readable early review of integrators for solar system dynamics is \cite{morbidelli02mis}, cf also \cite{morbidelli02mcm}}.

One of the first symplectic integrations of the solar system was done in \cite{sussman92ces} where it was confirmed that the solar system has a positive Lyapunov exponent, and hence exhibits chaotic behaviour cf \cite{laskar03css}.

More recently these methods have been improved and extended \cite{mclachlan95cmp,duncan98mts,laskar11nos,blanes15nfs}. Several symplectic integrators of high order were tested in \cite{farres13hps}, in order to determine the best splitting scheme for long-term studies of the solar system.

These various methods have resulted in the fact that numerical algorithms for solar system dynamics are now so accurate that they can be used to define the geologic time scales in terms of the initial conditions and parameters of solar system models (or vice versa).

\subsection{Symmetric Zassenhaus splitting and the equations of quantum mechanics} 

Equations of quantum mechanics in the semiclassical regime represent a double challenge of structure conservation and high oscillation.  A good starting point is the linear Schr\"odinger equation
\begin{equation}
  \label{LSE}
  \frac{\partial u}{\partial t}=\ii\varepsilon \frac{\partial^2 u}{\partial x^2}-\ii\varepsilon^{-1} V(x)u
\end{equation}
(for simplicity we restrict our discourse to a single space dimension), given in $[-1,1]$ with periodic boundary conditions. Here $V$ is the potential energy of a quantum system, $|u(x,t)|^2$ is a position density of a particle and $0<\varepsilon\ll1$ represents the difference in mass between an electron and  nuclei.  It is imperative to preserve the unitarity of the solution operator (otherwise $|u(\,\cdot\,,t)|^2$ is no longer a probability function), but also deal with oscillation at a frequency of $\O{\varepsilon^{-1}}$. A conventional approach advances the solution using  a palindromic splitting \R{palindromic}, but this is suboptimal for a number of reasons. Firstly, the number of splittings increases exponentially with order. Secondly, error constants are exceedingly large. Thirdly, quantifying the quality of approximation in terms of the step-size $h$ is misleading, because there are three small quantities at play: the step size $h$, $N^{-1}$ where $N$ is the number of degrees of freedom in space discretisation (typically either a spectral method or spectral collocation) and, finally, $\varepsilon>0$ which, originating in physics rather than being a numerical artefact, is the most important. We henceforth let $N=\O{\varepsilon^{-1}}$ (to resolve the high-frequency oscillations) and $h=\O{\varepsilon^\sigma}$ for some $\sigma>0$ -- obviously, the smaller $\sigma$, the larger the time step. 

\citeasnoun{bader14eas} have recently proposed an alternative approach to the splitting of \R{LSE}, of the form
\begin{equation}
  \label{Zassenhaus}
  \ee^{\ii h (\varepsilon \partial_x^2-\varepsilon^{-1} V)}\approx \ee^{\mathcal{R}_0}\ee^{\mathcal{R}_1} \cdots \ee^{\mathcal{R}_s}\ee^{\mathcal{T}_{s+1}}\ee^{\mathcal{R}_s} \cdots \ee^{\mathcal{R}_1} \ee^{\mathcal{R}_0}
\end{equation}
such that $\mathcal{R}_k=\O{\varepsilon^{\alpha_k}}$, $\mathcal{T}_{s+1}=\O{\varepsilon^{\alpha_{s+1}}}$, where $\alpha_0\leq \alpha_1<\alpha_2<\alpha_3<\cdots$ -- the {\em symmetric Zassenhaus splitting.\/} Here $\partial_x=\partial / \partial x$.

The splitting \R{Zassenhaus} is derived at the level of differential operators (i.e., prior to space discretisation), applying the symmetric Baker--Campbell--Hausdorff formula to elements in the free Lie algebra spanned by $\partial_x^2$ and $V$. For $\sigma=1$, for example, this yields
\begin{Eqnarray*}
    \mathcal{R}_0&=&-\Frac12\tau\varepsilon^{-1}V=\O{1},\\
    \mathcal{R}_1&=&\Frac12\tau\varepsilon \partial_x^2=\O{1},\\
    \mathcal{R}_2&=&\Frac{1}{24}\tau^3\varepsilon^{-1}(\partial_xV)^2+\Frac{1}{12} \tau^3\varepsilon \{(\partial_x^2V)\partial_x^2+\partial_x^2[(\partial_x^2V)\,\cdot\,]\}=\O{\varepsilon^2},\\
    \mathcal{R}_3&=&-\Frac{1}{120}\tau^5\varepsilon^{-1}(\partial_x^2 V)(\partial_xV)^2 -\Frac{1}{24}\tau^3\varepsilon (\partial_x^4V) +\Frac{1}{240}\tau^5\varepsilon \left(7\{(\partial_x^2V)^2\partial_x^2 \right.\\
    &&\mbox{}+\partial_x^2[(\partial_x^2V)^2\,\cdot\,] +\{(\partial_x^3V)(\partial_xV)\partial_x^2\left.\mbox{}+\partial_x^2[(\partial_x^3V)(\partial_xV)\,\cdot\,]\}\right)\\
    && +\Frac{1}{120}\tau^5\varepsilon^{-3} \{(\partial_x^4V)\partial_x^4+\partial_x^4[(\partial_x^4V)\,\cdot\,]\}=\O{\varepsilon^4},
\end{Eqnarray*}
where $\tau=\ii h$. Note that all the commutators, ubiquitous in the BCH formula, have disappeared: in general, the commutators in this free Lie algebra can be replaced by linear combinations of derivatives, with the remarkable property of {\em height reduction:\/} each commutator `kills' one derivative, e.g.
\begin{displaymath}
  [V,\partial_x^2]=-(\partial^2_x V)-2(\partial_xV)\partial_x,\qquad [[V,\partial_x^2],\partial_x^2]=(\partial_x^4V)+4(\partial_x^3V)\partial_x+4(\partial_x^2V)\partial_x^2.
\end{displaymath}

Once we discretise with spectral collocation, $\mathcal{R}_0$ becomes a diagonal matrix and its exponential is trivial, while $\ee^{\mathcal{R}_1}\MM{v}$ can be computed in two FFTs for any vector $\MM{v}$ because $\mathcal{R}_1$ is a Toeplitz circulant. Neither $\mathcal{R}_2$ nor $\mathcal{R}_3$ possess useful structure, except that they are {\em small!\/} Therefore we can approximate $\ee^{\mathcal{R}_k}\MM{v}$ using the Krylov--Arnoldi process in just 3 and 2 iterations for $k=2$ and $k=3$, respectively, to attain an error of $\O{\varepsilon^6}$ \cite{bader14eas}.

All this has been generalised to time-dependent potentials and is applicable to a wider range of quantum mechanics equations in the semiclassical regime.

\setcounter{equation}{0}
\setcounter{figure}{0}
\section{Beyond GNI}

Ideas in one area of mathematical endeavour often inspire work in another area. This is true not just because new mathematical research equips us with a range of innovative tools but because it provides insight that casts new light not just on the subject in question but elsewhere in the mathematical universe. GNI has thus contributed not just toward its own goal, better understanding of structure-preserving discretisation methods for differential equations, but in other, often unexpected, directions.

\subsection{GNI meets abstract algebra}

The traditional treatment of discretisation methods for differential equations was wholly analytic, using tools of functional analysis and approximation theory. (Lately, also tools from algebraic topology.) GNI has added an emphasis on geometry and this leads in a natural manner into concepts and tools from abstract algebra. As often in such mathematical dialogues, while GNI borrowed much of its conceptual background from abstract algebra, it has also contributed to the latter, not just with new applications but also new ideas.
\begin{itemize}
\item {\em B-series and beyond.\/} Consider numerical integration methods that associate to each vector field $\MM{f}$ a map $\MM{\psi}_h(\MM{f})$. A method $\MM{\psi}_h$ is called $g$-covariant\footnote{Also called equivariant.} if the following diagram commutes,
\begin{center}
  \begin{picture}(250,135)
     \thicklines
     \put (-20,0) {$\tilde{\MM{x}}=\MM{\psi}_h(\MM{f})(\MM{x})$} 
     \put (50,2) {\vector(1,0){140}}
     \put (200,0) {$\tilde{\MM{y}}=\MM{\psi}_h(\tilde{\MM{f}})(\MM{y})$}
     \put (15,110) {\vector(0,-1){95}}
     \put (225,110) {\vector(0,-1){95}}
     \put (-5,118) {$\dot{\MM{x}}=\MM{f}(\MM{x})$}
     \put (47,120) {\vector(1,0){148}}
     \put (210,118) {$\dot{\MM{y}}=\tilde{\MM{f}}(\MM{y})$}
     \put (100,8) {$\MM{x}=\MM{g}(\MM{y})$}
     \put (100,126) {$\MM{x}=\MM{g}(\MM{y})$}
  \end{picture}
\end{center}
It follows that if $g$ is a symmetry of the vector field $f$ and $\psi$ is $g$-covariant, then $\psi$ preserves the symmetry $g$. It seems that this concept of covariance for integration methods was first introduced in \cite{mclachlan95cps} and \cite{mclachlan98nit}.

It is not hard to check that all B-series methods are covariant with respect to the group of affine transformations.  A natural question to ask then, was ``are B-series methods the only numerical integration methods that preserve the affine group?". This question was open for many years, until it was answered in the negative  by \cite{munthekaas15abs}, who introduced a more general class of integration methods dubbed ``aromatic Butcher series", and showed that (under mild assumptions) this is the most general class of methods preserving affine covariance. Expansions of methods in this new class contain both rooted trees (as in B-series), as well as products of rooted trees and so-called $k$-loops \cite{iserles07bmc}. 

Whereas it may  be said that to date the importance of aromatic B-series has been at the formal rather than at the constructive level, these methods may hold the promise of the construction of affine-covariant volume-preserving integrators.

\item {\em Word expansions.\/} Classical B-series can be significantly generalised by expanding in {\em word series\/}   \cite{murua15wsd}. This introduced an overarching framework for Taylor expansions, Fourier expansions, modulated Fourier expansions and splitting methods. We consider an ODE of the form
\begin{equation}
  \label{word_series}
  \dot{\MM{x}}=\sum_{a\in\mathcal{A}} \lambda_a(t) \MM{f}_a(\MM{x}),\qquad \MM{x}(0)=\MM{x}_0,
\end{equation}
where $\mathcal{A}$ is a given {\em alphabet.\/} The solution of \R{word_series} can be formally expanded in the form
\begin{displaymath}
  \MM{x}(t)=\sum_{n=0}^\infty \sum_{\Mm{w}\in\mathcal{W}_n} \alpha_{\Mm{w}}(t) f_{\Mm{w}}(\MM{x}_0),
\end{displaymath}
where $\mathcal{W}_n$ is the set of all words with $n$ letters from $\mathcal{A}$. The coefficients $\alpha_{\Mm{w}}$ and functions $\MM{f}_{\Mm{w}}$ can be obtained recursively from the $\lambda_a$s and $\MM{f}_a$s in a manner similar to B-series. Needless to say, exactly like with B-series, word series can be interpreted using an algebra over rooted trees.

The concept of word series is fairly new in numerical mathematics but it exhibits an early promise to provide a powerful algebraic tool for the analysis of dynamical systems and their discretisation.

\item {\em Extension of Magnus expansions.\/} Let $\mathcal{W}$ be a {\em Rota--Baxter algebra,\/} a commutative unital algebra equipped with a linear map $R$ such that
\begin{displaymath}
  R(x)R(y)=R(R(x)y+xR(y)+\theta xy),\qquad x,y\in\mathcal{W},
\end{displaymath}
where $\theta$ is a parameter. The inverse $\partial$ of $R$ obeys
\begin{displaymath}
  \partial(xy)=\partial(x)y+x\partial(y)+\theta\partial(x)\partial(y)
\end{displaymath}
and is hence a generalisation of a derivation operator: a neat example, with clear numerical implications, is the backward difference $\partial(x)=[x(t)-x(t-\theta)]/\theta$. \citeasnoun{ebrahimifard09amf} generalised Magnus expansions to this and similar settings, e.g.\ dendriform algebras. Their work uses the approach in \cite{iserles99sld}, representing individual `Magnus terms' as rooted trees, but generalises it a great deal.

\item {\em The algebra of the Zassenhaus splitting.\/} The success of the Zassenhaus splitting \R{Zassenhaus} rests upon two features. Firstly, the replacement of commutators by simpler, more tractable expressions and, secondly, height reduction of derivatives under commutation. \citeasnoun{singh15ath} has derived an algebraic structure $\GG{J}$ which, encoding these two features, allows for a far-reaching generalisation of the Zassenhaus framework. The elements of $\GG{J}$ are operators of the form $\langle f\rangle_k =f\circ\,\partial_x^k+\partial_x^k\circ f$, where $k\in\BB{Z}_+$ and $f$ resides in a suitable function space. $\GG{J}$ can be endowed with a Lie-algebraic structure and, while bearing similarities with the Weyl algebra and the Heisenberg group, is a new and intriguing algebraic concept.
\end{itemize}

\subsection{Highly oscillatory quadrature} 

Magnus expansions \R{Magnus} are particularly effective when the matrix $A(t)$ oscillates rapidly. This might seem paradoxical -- we are all conditioned to expect high oscillation to be `difficult' -- but actually makes a great deal of sense. Standard numerical methods are based on Taylor expansions, hence on {\em differentiation,\/} and their error typically scales as a high derivative of the solution. Once a function oscillates rapidly, differentiation roughly corresponds to multiplying amplitude by frequency, high derivatives become large and so does the error. However, the Magnus expansion does not differentiate, it {\em integrates!\/} This has an opposite effect: the more we integrate, the smaller the amplitude and the series \R{Magnus} converges more rapidly. Indeed, often it pays to render a linear system highly oscillatory by a change of variables, in a manner described in  \cite{iserles02ged}, and then solve it considerably faster and cheaper. Yet, once we contemplate the discretisation of \R{Magnus} for a highly oscillatory matrix function $A(t)$, we soon come up another problem, usually considered difficult, if not insurmountable: computing multivariate integrals of highly oscillatory functions.

In a long list of  methods for highly oscillatory quadrature (HOQ) {\em circa\/} 2002, ranging from the useless to the dubious, one method stood out: \citeasnoun{levin82pco} proposed to calculate univariate integrals by converting the problem to an ODE and using collocation. This was the only effective method around, yet incompletely understood. 

The demands of GNI gave the initial spur to the emergence in the last ten years to a broad swath of new methods for HOQ: Filon-type methods, which replace the {\em non-oscillatory\/} portion of the integrand by an interpolating polynomial  \cite{iserles05eqh}, improved Levin-type methods \cite{olver06qmh} and the method of numerical stationary phase of \citeasnoun{huybrechs06eho}. The common characteristic of all these methods is that they are based on asymptotic expansions. This means that high oscillation is no longer the enemy -- indeed, the faster the oscillation, the smaller the error!

Highly oscillatory integrals occur in numerous applications, from electromagnetic and acoustic scattering to fluid dynamics, quantum mechanics and beyond. Their role in GNI is minor. However, their modern numerical theory was originally motivated by a problem in GNI. This is typical to how scholarship progresses and it is only natural that HOQ has severed its GNI moorings and has become an independent area on its own.

\subsection{Structured linear algebra}

GNI computations often lead to specialised problems in numerical linear algebra. However, structure preservation has wider impact in linear algebraic computations. Often a matrix in an algebraic problem belongs to an algebraic structure, e.g.\ a  specific Lie algebra or a symmetric space, and it is important to retain this in computation -- the sobriquet ``Geometric Numerical Algebra'' might be appropriate! Moreover, as in GNI so in GNA, respecting structure often leads to better, more accurate and cheaper numerical methods. Finally, structured algebraic  computation is often critical to GNI computations.
\begin{itemize}
\item Matrix factorization is the lifeblood of numerical algebra, the basis of the most effective algorithms for the solution of linear systems, computation of eigenvalues and solution of least-squares problems. A major question in GNA is ``Suppose that $A\in\mathcal{A}$, where $\mathcal{A}$ is a set of matrices of given structure. Given a factorization  $A=BC$ according to some set of rules, what can we say about the structure of $B$ or $C$?''. \citeasnoun{mackey05sfs} addressed three such `factorization rules': the {\em matrix square root,\/} $B=C$, the {\em matrix sign,\/} where the elements of $B$ are $\pm1$, and the {\em polar decomposition,\/} with unitary $B$ and positive semidefinite $C$. They focussed on sets $\mathcal{A}$ generated by a sesquilinear form $\langle\,\cdot\,,\,\cdot\,\rangle$. Such sets conveniently fit into two classes:
\begin{enumerate}
\item[(a)] Automorphisms $G$, such that $\langle G\MM{x},G\MM{y}\rangle=\langle\MM{x},\MM{y}\rangle$, generate a {\em Lie group;\/}
\item[(b)] Self-adjoint matrices $S$, such that $\langle S\MM{x},\MM{y}\rangle=\langle \MM{x},S\MM{y}\rangle$, generate a {\em Jordan algebra;\/} and
\item[(c)] Skew-adjoint matrices $H$ such that $\langle H\MM{x},\MM{y}\rangle=-\langle\MM{x},H\MM{y}\rangle$, generate a {\em Lie algebra.\/}
\end{enumerate}
It is natural to expect that conservation of structure under factorization would depend on the nature of the underlying inner product. The surprising outcome of \cite{mackey05sfs}   is that, for current purposes, it is sufficient to split sesquilinear forms into just two classes, unitary and orthosymmetric, each exhibiting similar behaviour.

\item Many algebraic eigenvalue problems are structured, the simplest example being that the eigenvalues of a symmetric matrix are real and of a skew-symmetric are pure imaginary: all standard methods for the computation of eigenvalues respect this. However, many other problems might have more elaborate structure, and this is the case in particular for nonlinear eigenvalue problems. An important example, with significant applications in mechanics, is 
\begin{equation}
  \label{QuadEig}
  (\lambda^2 M+\lambda G+K)\MM{x}=\MM{0},
\end{equation}
where both $M$ and $K$ are symmetric, while $G$ is skew symmetric. The eigenvalues $\lambda$ of \R{QuadEig} exhibit {\em Hamiltonian\/} pattern: if $\lambda$ is in the spectrum then so are $-\lambda,\bar{\lambda}$ and $-\bar{\lambda}$.\footnote{To connect this to the GNI narrative, such a pattern is displayed by matrices in the {\em symplectic Lie algebra\/} $\Gg{sp}(2n)$.} As often in numerical algebra, \R{QuadEig} is particularly relevant when the underlying matrices are large and sparse.

Numerical experiments demonstrate that standard methods for the computation of a quadratic eigenvalue problems may fail to retain the Hamiltonian structure of the spectrum but this can be obtained by bespoke algorithms, using a symplectic version of the familiar Lanczos algorithm, cf.\  \cite{benner07slc}. 

This is just one example of the growing field of structured eigenvalue and inverse eigenvalue problems.

\item The exponential from an algebra to a group: Recall Lie-group methods from Section~2.2: a critical step, e.g.\ in the RKMK methods, is the exponential map from a Lie algebra to a Lie group. Numerical analysis knows numerous effective ways to approximate the matrix exponential \cite{moler03ndw}, yet most of them fail to map a matrix from a Lie algebra to a Lie group! There is little point to expand intellectual and computational effort to preserve structure, only to abandon the latter in the ultimate step, and this explains the interest in the computation of the matrix exponential which is assured to map $A$ in a Lie algebra to an element in the corresponding Lie group. 

While early methods have used structure constants and, for maximal sparsity, Lie-algebraic bases given by space-root decomposition \cite{celledoni01mam}, the latest generation of algorithms is based upon {\em generalised polar decomposition} \cite{munthekaas01gpd}. 

\end{itemize}

\section*{Acknowledgments} 

This work has been supported by the Australian Research Council. The authors are grateful to David McLaren for assistance during the preparation of this paper, as well as to Philipp Bader,  Robert McLachlan and Marcus Webb, whose  comments helped to improve this paper.

\bibliographystyle{agsm}
\bibliography{refs}

\end{document}